\documentclass[reqno,11pt]{amsart}
\usepackage{latexsym,amssymb,amsfonts,amsmath,amsthm}
\usepackage{mathrsfs}
\usepackage[usenames]{color}
\usepackage{eucal}
\usepackage{graphicx}

\theoremstyle{definition}

\theoremstyle{remark}

\numberwithin{equation}{section}
\newcommand{\lms}{\mathop{\overline{\lim}}\limits}

\newcommand{\lmi}{\mathop{\underline{\lim}}\limits}
\def\v{\varepsilon}

\def\d{\delta}
\def\a{\alpha}
\def\b{\beta}
\def\g{\gamma}
\def\l{\lambda}

\def\n{\nu}

\def\x{\xi}

\def\e{\eta}

\def\C{\mathbb C}
\def\R{\mathbb R}

\def\V{\mathbb V}
\def\B{\mathbb B}
\def\F{\mathbb F}
\def\U{\mathbb U}

\def\r{\rho}
\def\x{\xi}

\begin{document}

 \title[Large Deviations for Processes on Half-Line] {Large Deviations for Processes on Half-Line: Random Walk and Compound Poisson.} %
 \author{F.C. Klebaner }
 \address{School of Mathematical Sciences,   Monash University,  Australia.}
  \email{fima.klebaner@monash.edu}
   \author{A.A. Mogulskii}
 \address{Sobolev Institute of Mathematics,
     Russia.}
  \email{mogul@math.nsc.ru}

\begin{abstract}
  We establish, under the Cramer exponential moment condition in a neighbourhood of zero, the Extended Large Deviation Principle for  the  Random Walk  and the Compound Poisson processes  in the  metric space $\V$ of functions of finite variation on $[0,\infty)$  with the modified Borovkov metric 
$\r(f,g)= \r_\B(\hat{f},\hat{g}) $, where $ \hat f(t)= f(t)/(1+t)$, $t\in \R$, and $\r_\B$ is the Borovkov metric. LDP in this space is ``more precise" than that with the usual metric of uniform convergence on compacts.
\end{abstract}

\thanks{This research was supported by the Russian Fund for Fundamental Research (projects number 14-01-0020), and the Australian Research Council Grant DP150103588.}
 \date{September 2016}%

\maketitle

\section{\bf Introduction}
The theory of Large Deviations for trajectories of processes seen as elements of the appropriate function space is well developed. However, for functions defined     on  infinite intervals, such as $\R^+$,   the typical metric used for LDP is that
of uniform convergence on compacts,  for example, for LDP for continuous processes the space of continuous functions $\C$ is used with the metric (e.g. \cite{Puh1}, \cite{Dem}, \cite{Kur})
\begin{equation}{\label{0.12}}
  \r^{(P)}(f,g):=\sum_{k=1}^\infty 2^{-k}\min\{\sup_{0\le t\le k}|f(t)-g(t)|, 1\}.
\end{equation}

Convergence  in metric $\r^{(P)}$
is equivalent to convergence in $\C[0,T]$ with uniform metric for any $T\geq 0$, e.g.    \cite{Puh2}. Hence a drawback of this  metric  is that it  is ``not sensitive" to the behaviour of functions at infinity. In \cite{BIB5} the LDP in the space $C$ with metric $
  \r (f,g)=\sup_{t\ge 0}\frac{|f(t)-g(t)|}{1+t},
$ on $\R^+$ is obtained for Diffusions and Random Walk, and is shown to be ``more precise" than in the space $(\C,\r^{(P)})$.
For discontinuous processes,
Dobrushin   and Pecherskij  \cite{BIB11}  give the LDP
for Compound Poisson processes
on the half-line using a metric based on the uniform metric.  Here we work under less stringent conditions, assuming exponential moments in the neighborhood of zero ($[{\bf C}_0]$), rather than on the whole line ($[{\bf C}_\infty]$), and   generalize their results   by using a different metric.
LDP for
Compound Poisson processes on $[0,1]$    in  the space of functions of bounded variation on $[0,1]$   is given in \cite{BIB10}  and recently  generalized in \cite{BIB4}   by using  Borovkov's metric $\r_\B$ instead of uniform.
Here we extend    results of \cite{BIB4} to the half-line yielding a  generalization of \cite{BIB11}. Since here we work under less stringent moment conditions $[{\bf C}_0]$, we also generalize results of \cite{BIB5}, that give
  the classical LDP  on the half-line  for Random Walk under $[{\bf C}_\infty]$. The proofs are different to those in \cite{BIB5}.

For $k\ge 0$, denote by $S_k=\sum_{i=1}^k\x_{(i)}$ the partial sums  of  i.i.d. r.v.'s $\{\x_{(i)}\}_{i\ge 1}$  distributed as
$\x$, and $S_0=0$.
Let $s=s(t)$ be the continuous piecewise linear function
on $[0,\infty)$ going
  through  the points
 $
  (0,S_0), (1,S_1),\cdots, (k,S_k),\cdots
 $
Define the process $s_n$ by
$$
   s_n=s_n(t):=\frac{1}{x}s(tn),~~~t\ge 0,~~~n=1,2,\cdots,
$$
where  $x=x(n)\sim n$ as $n\to \infty$.


Similarly the process   $\x_T$  for a real   $T\ge 1$     is defined.   Consider a Compound Poisson process
$$
  \x=\x(t),~~~t\ge 0,
$$
and let
\begin{equation}{\label{1.2}}
  \x_T= \x_T(t):= \frac{1}{x}\x(tT),~~~t\ge 0,~~~T\ge 1,
\end{equation}
where   $x=x(T)\sim T$ as $T\to \infty$.

Two families $s_n$ and $\x_T$ (processes $ s(t)$ and $\x (t)$) have much in common.     $\x(t)$
has independent increments, and so does  $s(t)$, when taken at integer times $t=n$.
If the r.v.  $\x$ in the definition of  $s(t)$  is taken to be $\x(1)$, then   $s(t)$ is   the linear interpolation of  $\x(t)$, going through the points
$$
  (0,\x(0)), (1,\x(1)), \cdots, (k,\x(k)), \cdots
$$
Therefore it is not surprising, that under mild assumptions, the families
$s_n$ and $\x_T$
satisfy  Large Deviation Principle with the common Rate Function, determined by the r.v.  $\x$.
We assume throughout that  r.v.
$\x$ in definition (\ref{1.1})  satisfies   Cramer condition $[{\bf C}_0]$.

$[{\bf C}_0]$. {\it For some $\d>0$}
$$
  {\bf E}e^{\d|\x|}<\infty.
$$

 We establish the  {\it Extended Large Deviation Principle} for   two families
\begin{equation}{\label{1.1}}
  \{s_n=s_n(t);~t\ge 0\}_{n\ge 1},\;\;\mbox{and\;\;}\{\x_T=\x_T(t);~t\ge 0\}_{T\ge 1}
\end{equation}
defined on the half-line, $t\in [0,\infty)$.

%
%

   The precise definition of ELDP is given in Section 3,  (see also \cite{BIB1} or \cite{BIB2}, ch. 4).
 ELDP      holds under less stringent requirements on the rate function than the classical LDP, in particular the space is not required to be complete,  and the rate function, while lower semi-continuous, is not required to be compact. However, if ELDP holds with a good rate function, then LDP follows.

  ELDP for    processes in (\ref{1.1}) defined on   $[0,1]$ under the  assumption $[{\bf C}_0]$   in the space of functions    on   $[0,1]$  was established earlier in
\cite{BIB3}  and
 \cite{BIB4} (see also \cite{BIB2}, ch. 4). The main contribution of this work is to extend the  results of
 \cite{BIB3} and
 \cite{BIB4} to processes defined
on   $[0,\infty)$ and establish ELDP in the space of functions defined on the half line.

An extension of   the classical LDP    to the half-line    was recently done in \cite{BIB5}, in particular for Random Walk $s_n$, but under a stronger Cramer  condition $[{\bf C}_\infty]$.

$[{\bf C}_\infty]$. {\it For all $\l\in \R$}
 $$
   {\bf E}e^{\l\x}<\infty.
 $$
Note  here that     $[{\bf C}_\infty]$ is a necessary and sufficient condition for the classical LDP for $s_n$ in the metric space $\C[0,1]$ with the uniform metric (by   Puhalski's Theorem LDP is equivalent to exponential tightness which is equivalent to  $[{\bf C}_\infty]$,     eg.  Lemma 4.4.5 in \cite{BIB2}).

The paper is organised as follows.   In Section 2 we introduce the   metric space   of functions    of bounded variation $\V$  defined on
$[0,\infty)$, with the metric $\r$,   based on the Borovkov's metric. Section  3 contains main definitions and  results. Sections 4-9 contain   proofs.

\section{\bf The space $(\V,\r)$.}
\setcounter{section}{2}
\setcounter{equation}{0}

We look at processes $s_n=s_n(t)$ and
$\x_T=\x_T(t)$  as random elements of the space  $\V$
of functions
 $f=f(t)$, defined for $t\in \R$,
   {\it having bounded variation on any interval $[0,T]$,
   without discontinuities   of the second kind, such that $f(t)=0$ for $t\le 0$;  at a point of discontinuity $t_0$  the function
  $f$   can take any value $f(t_0)$, in the interval
  $[f(t_0-,~f(t_0+))]$}.
Define  $\x_T(t)$ for $t\le 0$ accordingly,
 $$
    \x_T=\x_T(t):= \begin{cases}
         0,~~~t\le 0;\\
       \frac{1}{x}\x(tT),~~~t\ge 0.
 \end{cases}
  $$
Similarly define $s_n(t)=0$ for $t\le 0$.
It is well known, (eg. \cite{BIB6}), that it is possible to define the process $\x=\x(t)$
 to be left-continuous, so that
 $$
  {\bf P} (\x_T\in \V)=1~~~\mbox{for all}~~~T>0.
 $$
For every $f\in \V$ consider its   {\it graph} $\Gamma_f$, a simply connected set in $\R^2$,  determined by its sections  $$
  \Gamma_f|_u:=\Gamma_f\cap\{(t,\a)\in \R^2:~t=u\}= [[(u,f(u- )),~(u,f(u+ ))]],
 $$
where $[[\a,~\b]]$ denotes the line connecting
 two points
 $\a,\b\in \R^2$. If $f $ is  continuous at $t=u$  then $$
    \Gamma_f|_u=(u,f(u))
 $$
 consists of the single point $(u,f(u))$. If $t=u$
is a point of discontinuity of  $f $, then the section of the graph at this point is a vertical line segment $$
    \Gamma_f|_u=\{(u,\a):~~\a\in [\a_-,\a_+]\},~~~\mbox{where}~~~\a_-:=f(u- ),
    ~~\a_+:=f(u+ ).
 $$

It is convenient to use the ``square'' norm in $\R^2$,  in which   the graphs of  $f\in \V$ lie,
 $[(t,\a)]:=\max\{|t|,|\a|\}$
 ($\v$-neighborhood   $(\g)_\v$  of any point $\g=(t,\a)\in \R^2$, in this norm is a square with the centre at
 $\g$, with sides parallel to the coordinates and length $2\v$). For a set  $A\subset \R^2$ denote by  $(A)_\v$
 its $\v$-neighborhood in this ``square'' norm.

We use the Borovkov's metric  in the space  $\V$, $\r_\B=\r_\B(f,g)$,  defined as follows
  $\r_\B(f,g)<\v$   {\it if and only if both relations hold}
$$
\Gamma_f\subset (\Gamma_g)_\v~~~\mbox{and}~~~\Gamma_g\subset (\Gamma_f)_\v.
$$

The metric $\r_\B$ was introduced by Borovkov A.A. in \cite{BIB7}
(for the space $\F$, wider than $\V$, see also \cite{BIB8}); the topology generated by
$\r_\B$, is the same as Skorohod $M_2$ topology, described in
\cite{BIB9}.  The metric    $\r_\B$ was effectively used in
\cite{BIB3}, \cite{BIB4}.
Note that $\r_\B$  is weaker than the uniform metric $\r_\U$, $\r_\U(f,g):=\sup_{t}|f(t)-g(t)|$, i.e.
\begin{equation}{\label{2.1}}
  \r_\B(f,g)\le \r_\U(f,g)~~~\mbox{for all}~~~f,g\in \V.
\end{equation}
Indeed, for any  $t\in \R$
the point  $(t,\a)$ in  $\Gamma_f$  is given by
$$
  (t,\a)=(t,pf(t- )+qf(t+ )),~~~\mbox{where}~~~p\ge 0,~~q\ge 0,~~p+q=1,
$$
moreover  $|f(t\pm  )-g(t\pm  )|\le \r_\U(f,g)$,
so that (\ref{2.1}) holds.

The main metric for our analysis   $\r=\r(f,g)$ is obtained from the Borovkov's metric
$\r_\B$ by weighing functions $f$ and $g$
$$
  \hat{f}(t):=\frac{f(t)}{1+|t|},~~~\hat{g}(t):=\frac{g(t)}{1+|t|},~~~t\in \R,
$$
$$\r(f,g):= \r_\B(\hat{f},\hat{g}).$$

 LDP  for the family $s_n $ on  the half-line under condition $[{\bf C}_\infty]$ is given in  \cite{BIB5}
 where another metric $\hat{\r}$ was used, obtained by weighing functions in the uniform metric
 \begin{equation}{\label{2.3}}
  \hat{\r}(f,g)=\r_\U(\hat{f},\hat{g}).
\end{equation}
 $\r$ is weaker than  $\hat{\r}$, since
by (\ref{2.1})
\begin{equation}{\label{2.2}}
  \r(f,g)\le \hat{\r}(f,g)~~~\mbox{for any}~~~f,g\in \V,
\end{equation}
Denote by $\V^0\subset\V$
the class of functions $f\in \V$,
such that
$$
  \lim_{t\to \infty}\frac{|f(t)|}{1+t}=\lim_{t\to \infty}|\hat{f}(t)|=0.
$$

To summarise,   the families of   processes $\{s_n\}$, $\{\x_T\}$ have trajectories
in $\V$. The distribution of $s_n$ ($\x_T$)
is determined by the norming sequence
$x=x(n)\sim n$  ($x=x(T)\sim T$) and the distribution of the rv.  $\x$,
that denotes the jump in the random walk (the increment of  $\x(t)$
on the unit interval  $\x(1)$).
   In this way   $\x$ denotes  two different random variables from the two families. The main moment conditions   and the rate function  are given  in terms of
$\x$ (Section 3).

Without loss of generality, by changing   the drift  if necessary, we can assume  ${\bf E}\x=0$. In this case
  the trajectories of   $s_n$ and $\x_T$,
belong to
$\V^0\subset\V$ with probability one. As it will be seen in
Lemma  3.1,
any  $f\in \V$ with  $J(f)<\infty$
belongs to $\V^0$.  However, the main Theorem 3.1 uses the space $\V$ (while we could have used
$\V^0$).

\section{\bf Statements of main results. The rate function $J=J(f)$.}
\setcounter{section}{3}
\setcounter{equation}{0}

Let $\x$  be a non-degenerate rv. with
${\bf E}\x=0$,  satisfying the Cramer condition $[{\bf C}_0]$.
Let  $\psi(\l)$ be the Laplace transform of
  $\x$
$$
  \psi(\l):={\bf E}e^{\l\x},~~~\l\in \R,
$$
and denote by $(\l_-,\l_+)$  the largest interval for which
$\psi(\l)$ is finite. Due to condition $[{\bf C}_0]$, this interval is not empty and contains the point
  $\l=0$.
Denote the Legendre transform of $ \log \psi(\l)$ (the deviation function of $\x$) by
$$
  \Lambda(\a):=\sup_\l\{\l\a-\log \psi(\l)\},~~~\a\in \R.
$$
The properties of $ \Lambda$ are well known,
it is non-negative, convex, lower semi-continuous, equals to zero at a single point
$a={\bf E}\x$ (in our case $a=0$),   (e.g. \cite{Dem}, \cite{BIB6} or
\cite{BIB2}, ch 2).

Recall the decomposition of any $f\in \V$ into absolutely continuous and singular components  $f_a$ and $f_s$,
$$
  f=f_a+f_s=f_a+f_{s+}-f_{s-},~~~f_a\in \V,~~~f_{s\pm}\in \V,
$$
 where  $f_{s+}$  and  $-f_{s-}$ are non-decreasing and non-increasing parts of $f_s$.
Using this representation define the   functional  (cf. \cite{BIB2}, ch.4)  for any
$U\in (0,\infty)$
\begin{equation}\label{J0U}
  J_0^U(f):=\int_0^U\Lambda(f'_a(t))dt+\l_+f_{s+}(U) +
  |\l_-|f_{s-}(U).
\end{equation}
It is clear that
$J_0^U(f)$ is non-decreasing in $U$, therefore there is a limit as $U\to\infty$, which defines  the rate function
$$
  J(f):=\lim_{U\to \infty}J_0^U(f),~~~f\in \V.
$$

The properties of
$J(f)$ are summarised in the following Lemma 3.1.

\bigskip

{\bf Lemma 3.1.} {\it

$I.$  $J(f)$ is lower semi-continuous in the space $(\V,\r)$:
\begin{equation}{\label{3.1}}
  \lmi_{\r(f_n,f)\to 0}J(f_n)\ge J(f).
\end{equation}

$II.$  For some $C<\infty$ and all $f\in \V$,  if $J(f)\le N$  then
$$
  |f(U)|\le C\sqrt{U}N,~~~U\ge 1.
$$

$III.$  For any $f\in \V$ there exists a sequence of absolutely continuous functions
$f_n\in \V$ such that
$$
  \r(f_n,f)\to 0,\;\;\mbox{and}\;\;J(f_n)\to J(f), \; \mbox{ as }\; n\to \infty.
$$
}

The proof of Lemma 3.1 is given in Section 5, and     now we turn to the main result.
Denote, as usual, for a measurable non empty set  $B\subset \V$
 $$
  J(B):=\inf_{f\in B}J(f),~~~J(\emptyset)=\infty.
$$
Denote  by $(B)$, $[B]$   the interior and the closure of $B$ respectively,
and  by $(B)_\v$ the
$\v$-neighbourhood of   $B$  with respect to our metric $\r$.
Finally, let
$$
  J(B+):=\lim_{\v\downarrow 0}J((B)_\v).
$$
Since for any $\v>0$, the following inclusions hold
$$
  (B)\subset B \subset [B]\subset (B)_\v,
$$
we have that
$$
  J((B))\ge J(B)\ge J([B])\ge J(B+).
$$

\bigskip


{\bf Theorem 3.1}. {\it I. The family $s_n$ satisfies Extended Large Deviation Principle in the space  $(\V,\r)$  with the rate function $J$, namely   for any measurable set
  $B\subset \V$
$$
\lms_{n\to \infty}\frac{1}{n}\log {\bf P}(s_n\in B)\le -J(B+),
$$

$$
\lmi_{n\to \infty}\frac{1}{n}\log{\bf P}(s_n\in B)\ge -J((B)).
$$

II.  The family $\x_T$ satisfies Extended Large Deviation Principle in the space  $(\V,\r)$  with the rate function $J$, namely  for any measurable set
  $B\subset \V$
\begin{equation}{\label{3.2}}
\lms_{T\to \infty}\frac{1}{T}\log{\bf P}(\x_T\in B)\le -J(B+),
\end{equation}
\begin{equation}{\label{3.3}}
\lmi_{T\to \infty}\frac{1}{T}\log{\bf P}(\x_T\in B)\ge -J((B)).
\end{equation}
}

The proof  of Theorem 3.1  is given in Section 4.

  Now let us compare this result  with previously known results.
Large Deviation Principles for Compound Poisson processes on $[0,1]$  were established
 in \cite{BIB10} and \cite{BIB4}
under  Cramer condition $[{\bf C}_0]$ in   the space  $\V[0,1]$ of functions of bounded variation on $[0,1]$  with different metrics.
In  \cite{BIB10}   uniform metric $\r_\U$ is used, while
in
\cite{BIB4}
Borovkov's metric $\r_\B$ is used.
The main result of
  \cite{BIB4}
 strengthens, in particular, the main result of  \cite{BIB10}
for an important class of boundary crossing sets $$
  B_c:=\{f\in \V:~\sup_{0\le t\le 1}f(t)\ge c\},~~~c>0.
$$
Indeed, since the sequence
$$
f_n(t):= \begin{cases}
~~0,~~~0\le t\le \frac{1}{2}-\frac{1}{n};\\
1,~~~                  \frac{1}{2}-\frac{1}{n}<t\le \frac{1}{2}\\
0,~~~                  \frac{1}{2}<t\le 1,
\end{cases}
$$
converges weakly to $f_0:=f_0(t)\equiv 0$,
 $f_0$ belongs to the closure  $[B_1]$ of $B_1$ in the topology of weak convergence.  $J(f_0)=0$ (when the mean of the underlying process is zero)
giving a non informative   trivial upper bound  in  \cite{BIB10},  different
to the lower bound.
It can be seen (Lemma 10.1 of Appendix)  that
$$
 J((B_1))=J(B_1+)=\inf_{0\le v\le 1}J(g_v),
$$
  where $g_v$ is   a continuous piece-wise linear function given by  $g_v(t)=\frac{t}{v}$ for $t\in [0,v]$, and
 $g_v'(t)=a$   with   $a={\bf E}\x(1)$  for $t\in (v,1]$.
Therefore ELDP  in \cite{BIB4}, Theorem 1.1 allows to obtain the ``correct'' logarithmic asymptotic for   probability of $B_1$.

The paper \cite{BIB11}  generalizes the LDP
for Compound Poisson processes in  \cite{BIB10}
from the interval   to the half-line using a metric based on the uniform metric.
Here we generalize  ELDP of \cite{BIB4}  from the interval   to the half-line  using   a metric based on the Borovkov's metric.
The above illustration of   LDP's  with different metrics  shows advantages of our   result  as compared to that in  \cite{BIB11}.

\bigskip

If  the underlying random variable $\x$
satisfies  a  stronger  Cramer  condition  $[{\bf C}_\infty]$ instead of  $[{\bf C}_0]$, then our result implies the classical LDP on the half-line.
 Indeed, in this case $\l_+=|\l_-|=\infty$, and the rate function becomes
$$
  J(f)=I(f):=
    \begin{cases}
        \int_0^\infty\Lambda(f'(t))dt,~~~\mbox{if}~~f~~\mbox{is absolutely continuous};\\
       \infty,~~~\mbox{otherwise}.
 \end{cases}
$$
Further, as shown in   \cite{BIB5}, the rate function $I(f)$ is a
``good'' rate function in the space
of continuous  functions  on the half-line $(\C,\hat{\r})$ with the metric
$\hat{\r}$   in (\ref{2.3}), constructed by using the uniform metric $\r_\U$. This means that for any $v\ge 0$ the set
$\{f\in \C:~I(f)\le v\}$ is a compact in $(\C,\hat{\r})$. Since $\r$ is weaker than $\hat{\r}$, see  (\ref{2.2}),
it is clear   that $I$ remains a ``good'' rate function in the space  $(\V,\r)$. If for any $v$ the set  $\{f\in \V:~I(f)\le v\}$ is  a  compact, then
$$
  I(B+)=I([B]),
$$
where    $[B]$ is the closure of  $B$  in $(\V,\r)$, this is shown in \cite{BIB1} (see also Lemma  4.1.1 in
    \cite{BIB2}).
Therefore   Theorem 3.1  implies LDP on half-line for $\x_T$ as well as recovers the LDP for $s_n$ in $(\C,\hat{\r})$,     given recently in \cite{BIB5}.

 \bigskip

{\bf Corollary 3.1}. {\it I. Let r.v. $\x$
satisfy $[{\bf C}_\infty]$, ${\bf E}\x=0$. Then the family  $s_n$ satisfies LDP in the space $(\V,\r)$ with rate function $I$: for any
measurable set $B\subset \V$
$$
\lms_{n\to \infty}\frac{1}{n}\log{\bf P}(s_n\in B)\le -I([B]),
$$

$$
\lmi_{n\to \infty}\frac{1}{n}\log{\bf P}(s_n\in B)\ge -I((B)).
$$

II. Let r.v. $\x$
satisfy $[{\bf C}_\infty]$, ${\bf E}\x=0$. Then the family  $\x_T$ satisfies LDP in the space $(\V,\r)$ with rate function $I$: for any
measurable set $B\subset \V$
$$
\lms_{T\to \infty}\frac{1}{T}\log{\bf P}(\x_T\in B)\le -I([B]),
$$

$$
\lmi_{T\to \infty}\frac{1}{T}\log{\bf P}(\x_T\in B)\ge -I((B)).
$$

}

\section{\bf Proof of Theorem 3.1. }   We prove only the second statement for Compound Poisson processes.
Part $I$, for Random Walk, has a similar proof,    being    simpler in places.  Since the proof uses results for ELDP on compacts,
one needs to replace  references to  results in  \cite{BIB4}, where
ELDP is established for
$\x_T$  on  $[0,1]$, by  reference to   \cite{BIB3} (\cite{BIB2}), where ELDP is established for $s_n$ on  $[0,1]$).

\bigskip

The proof is based on Lemmas 4.1, 4.2 and 4.3  below.

\bigskip

{\bf Lemma 4.1.} {\it For any $\v>0$ and $f\in \V$

\begin{equation}{\label{4.1}}
\lmi_{T\to \infty}\frac{1}{T}\log {\bf P}(\x_T\in (f)_\v)\ge -J(f).
\end{equation}
}

 {\bf Lemma 4.2.} {\it For any $\v\in (0,~\frac{1}{10})$ and $f\in \V$

\begin{equation}{\label{4.2}}
\lms_{T\to \infty}\frac{1}{T}\log{\bf P}(\x_T\in (f)_\v)\le -J((f)_{59\v}).
\end{equation}
}

{\bf Lemma 4.3.} {\it For any $\v>0$ and $N<\infty$
there are $M<\infty$ and a collection of functions $\{g_1,\cdots,g_M\}$, $g_i\in\V$ so that
\begin{equation}{\label{4.3}}
\lms_{T\to \infty}\frac{1}{T}\log {\bf P}(\x_T\not \in \cup_{k=1}^M(g_k)_\v)\le -N.
\end{equation}
}

We prove these Lemmas later, and now give the proof  of the theorem.

\bigskip

\begin{proof} of Theorem 3.1.

 $II$.
$(i)$. Upper bound. For any $\v>0$  and $N<\infty$
by Lemma 4.3 there are functions  $\{g_1,\cdots, g_M\}$ in   $\V$,
such that (\ref{4.3}) holds.
Denote  ${\mathcal M}:=\{i\in \{1,\cdots,M\}:~B\cap (g_i)_\v\not = \emptyset\}$.
The following bound clearly  holds
\begin{equation}{\label{4.4}}
{\bf P}(\x_T\in B)\le {\bf P}(\x_T\not\in \cup_{i=1}^M(g_i)_\v)+
\sum_{i\in {\mathcal M}}{\bf P}(\x_T\in (g_i)_\v).
\end{equation}
Further, by Lemma  4.2, for $i\in {\mathcal M}$
\begin{equation}{\label{4.5}}
\lms_{T\to \infty}\frac{1}{T}\log {\bf P}(\x_T\in (g_i)_\v)\le -J((g_i)_{59\v}).
\end{equation}
By using  (\ref{4.3}), (\ref{4.5}),
 we obtain from (\ref{4.4}) the bound
\begin{equation}{\label{4.6}}
\lms_{T\to \infty}\frac{1}{T}\log {\bf P}(\x_T\in B)\le
-\min\{N, \min_{i\in {\mathcal M}}J((g_i)_{59\v}) \}.
\end{equation}
Since for any  $i\in {\mathcal M}$
 there is $f_i\in B$   such that $\r(f_i,g_i)<\v$,
we have
$$
  (g_i)_{59\v}\subset (f_i)_{60\v}\subset(B)_{60\v}.
$$
Therefore
$$
     \min_{i\in {\mathcal M}}J((g_i)_{59\v})\ge J((B)_{60\v}),
$$
and now it follows from   (\ref{4.6})
\begin{equation}{\label{4.7}}
\lms_{T\to \infty}\frac{1}{T}\log {\bf P}(\x_T\in B)\le
-\min\{N, J((B)_{60\v})\}.
\end{equation}
Taking   $N\to \infty$ and then $\v\to 0$ in (\ref{4.7}),
we obtain (\ref{3.2}).

$(ii)$. Lower bound (\ref{3.3}) follows directly from Lemma 4.1. Theorem 3.1 is proved.
\end{proof}

\section{\bf Proof  of properties of the rate function  Lemma  3.1.}
\begin{proof}
$I$.  Lower-semicontinuity of $J$. From the definition of $J(f)$ it follows that for any $N<\infty$ and $\d>0$
there is  $U=U_{N,\d}<\infty$
so that
$$
  J_0^U(f)\ge \min\{J(f)-\d,N\}.
$$
($N$ is used in case $J(f)=\infty$).
We can see by using   lower semi-continuity of  $J_0^U(f)$
in space $(\V,~\r)$ that for  $V>U$ and  $\r(f_n,f)\to 0$
\begin{equation}{\label{5.1}}
\lmi_{n\to \infty}J_0^V(f_n) \ge J_0^U(f).
\end{equation}
Now we have
$$
 \lmi_{n\to \infty}J(f_n) \ge \lmi_{n\to \infty}J_0^V(f_n) \ge J_0^U(f)\ge
  \min\{J(f)-\d,N\}.
$$
Since $N<\infty$ and $\d>0$ are arbitrary,
(\ref{3.1}) follows.

It remains to show  (\ref{5.1}).  The proof is similar to that of   of Theorem  4.2.2  part
$(ii)$ in  \cite{BIB2} for lower semi-continuity of $J$ when $f$ belongs to a space of functions defined on $[0,1]$. Since we work on $\R^+$   there are differences, and we give it here.
Denote by ${\bf t}_K=(t_0,\cdots,t_K)$, $0=t_0<t_1<\cdots<t_K=U$
a partition of  $[0,U]$
into  $K$ parts. For a function $g\in \V$,
  $g^{{\bf t}_K}$
denotes the continuous pice-wise linear function on $[0,U]$
going through the points
$$
  (t_k,~g(t_k)),~~~k=0,\cdots,K.
$$
Then, by definition
$$
  J_0^U(f^{{\bf t}_K}) = \int_0^U\Lambda((f^{{\bf t}_K})'(t))dt.
$$
 Theorem 4.2.1 of  \cite{BIB2} states that
\begin{equation}\label{Thm421}
J(f)=\sup_{{\bf t}_K} I(f^{{\bf t}_K}),
\end{equation}
where $\sup$ is over all partitions of $[0,U]$.
 Therefore
 for any $N<\infty$, $\d>0$
there is a partition
 ${\bf t}_K$ such that
\begin{equation}{\label{5.2}}
J_0^U(f^{{\bf t}_K})\ge \min\{J_0^U(f)-\d,N\}.
\end{equation}
For this partition we have,  due to $\r(f_n,f)\to 0$
as $n\to \infty$,   that for any  $k\in \{0,1,\cdots,K\}$
there is $(t_k^{(n)}, \hat{\a}_k^{(n)})\in \Gamma_{\hat{f_n}} $  such that
\begin{equation}{\label{5.3}}
      |t_k-t_k^{(n)}|\to 0,~~~|\frac{f(t_k)}{1+t_k}-\hat{\a}_k^{(n)}|\to 0~~~\mbox{as}~~~
      n\to \infty.
\end{equation}
Construct now the function  $g_n$ from $f_n$ by replacing  its values at $t_k^{(n)}$,  $f(t_k^{(n)})$ by \\
$g(t_k^{(n)}):=(1+t_k^{(n)})\hat{\a}_k^{(n)}$.
Clearly, the graph of  $f_n$  does not change, therefore by
\eqref{Thm421}
we have for  $V\ge V_N:=t_K^{(n)}$
\begin{equation}{\label{5.4}}
  J_0^V(f_n)=J_0^V(g_n)\ge
    J_0^{V_n}(g_n)\ge J_0^{V_n} (g_n^{{\bf t}^{(n)}_K}).
\end{equation}
From (\ref{5.3}) it follows that
$$
  \max_{0\le k\le t_K}\{|t_k-t_k^{(n)}|\} \to 0,~~~
  \max_{0\le k\le t_K}\{|f(t_k)-g_n(t_k^{(n)})|\}\to 0~~~\mbox{as}~~~n\to \infty.
$$
Hence by lower semi-continuity of $\Lambda(\a)$
we have
\begin{equation}{\label{5.5}}
  \lmi_{n\to \infty} J_0^{V_n} (g_n^{{\bf t}^{(n)}_K})\ge
  J_0^U (f^{{\bf t}_K}).
\end{equation}
Now      (\ref{5.4}), (\ref{5.5}), (\ref{5.2})  imply
\begin{equation}{\label{5.6}}
  \lmi_{n\to \infty} J_0^V (f_n)\ge \min\{J_0^U(f)-\d,N\}.
\end{equation}
Thus  (\ref{5.1}) is established.

\bigskip

$II$.   For a fixed  $U\ge 1$  consider the set $$
  B_U:=\{g\in \V: g(U)=f(U)\}.
$$
 $J_0^U(g)$  achieves its minimum over
$g\in B_U$   on the function $g_0(t)=f(U)\frac{t}{U}$
for $0\le t\le U$, by equation  (\ref{10.3}),   Appendix, Lemma 10.1.
Therefore
$$
  N\ge J(f) \ge J_0^U(f)\ge J_0^U(g_0)=U\Lambda(\frac{f(U)}{U}).
$$
In view of $[{\bf C}_0]$ and  ${\bf E}\x=0$,
for some $c>0$
$$
  \Lambda(\a)\ge c\min\{\a^2,|\a|\},
$$
hence
$$
  N\ge c\min\{\frac{f^2(U)}{U},|f(U)|\}.
$$
It now follows that if
$$
  c\frac{|f(U)|^2}{U}\le N,~~~|f(U)|\le \frac{1}{\sqrt{c}}\sqrt{UN};
$$
if $\frac{|f|^2}{U}>|f(U)|$, then
$$
  c|f(U)|\le N,~~~|f(U)|\le \frac{1}{c}N.
$$
Hence for $T\ge 1$  we have
$$
   |f(U)|\le \max\{\frac{1}{\sqrt{c}}\sqrt{UN},\frac{1}{c}N\}\le C\sqrt{U}N,
$$
where $C:=\frac{1}{\sqrt{c}}+\frac{1}{c}$.
Statement  $II$ is proved.

$III$. If $J(f)=\infty$, then since any $f\in \V$ is a limit of absolutely continuous    $f_n\in \V$,
$\r(f_n,f)\to 0$ as $n\to \infty$, and using lower semi-continuity established in
  $I$, $$
  \lim_{n\to \infty}J(f_n)=\infty=J(f).
$$
Take $J(f)<\infty$.
Using the already proven results $I$ and $II$  we can show (following the proof of Theorem 4.2.1 in
  \cite{BIB2})
that for any $\d>0$
there is a sequence of piecewise linear $f_n\in \V$
going through the points
$$
  (0,0), (t^{(n)}_1,f(t^{(n)}_1)),\cdots,  (t^{(n)}_{K_n},f(t^{(n)}_{K_n})),~~~
  f_n(t)=f(t^{(n)}_{K_n})~~~\mbox{for}~~~t\ge t^{(n)}_{K_n},
$$
such that
$$
  \r(f_n,f) \le \frac{1}{n},~~~J(f_n)=J_0^{t^{(n)}_{K_n}}(f_n)\le J(f)+\frac{1}{n}.
$$
Therefore for this sequence
$$
  \lms_{n\to \infty} J(f_n)\le J(f),
$$
which together with (3.1) gives
$$
  \lim_{n\to \infty} J(f_n) = J(f).
$$
Lemma 3.1 is proved.

\end{proof}

\section{\bf  Auxillary statements.} To proceed we need the following results.
Denote for $U\ge 2$, $V\ge 2$
$$
  A(U,\v):=\{g\in \V:~\sup_{t\ge U}|\hat{g}(t)|\le \v\},
  $$
  $$
  B(2U,V):=\{g\in \V:~\sup_{t\le 2U}|\hat{g}(t)|\le V\}.
$$

{\bf Lemma 6.1.}
\begin{enumerate}
\item
 For any  $N<\infty$ and $\v>0$ there is  $U=U_{N,\v}<\infty$
such that
\begin{equation}{\label{6.1}}
\lms_{T\to \infty}\frac{1}{T}\log
{\bf P}(\x_T\not \in A(U,\v))\le -N.
\end{equation}
\item
  For any
 $N<\infty$ and $U<\infty$ there is $V=V_{N,U}<\infty$
such that
\begin{equation}{\label{6.2}}
\lms_{T\to \infty}\frac{1}{T}\log
{\bf P}(\x_T\not \in B(2U,V))\le -N.
\end{equation}
\end{enumerate}

\begin{proof} of Lemma 6.1.

$(1)$  Since
$$
  \sup_{t\ge U}\frac{1}{1+t}|\x_T(t)|\le
  \sup_{tT\ge UT}\frac{1}{tT}|\x(tT)| =
  \sup_{v\ge UT}\frac{1}{v}|\x(v)|,
$$
we have
\begin{equation}{\label{6.3}}
\{\x_T\not \in A_T(U,\v)\}=
\{\sup_{t\ge U}\frac{1}{1+t}|\x_T(t)| > \v\}\subset
\cup_{k\ge [UT]}B_k(\v),
\end{equation}
where
$$
  B_k(\v):=\{\sup_{k\le v\le k+1}\frac{1}{v}|\x(v)|\ge \v\}.
$$
Since
$$
   B_k(\v)=B^+_k(\v) \cup B^-_k(\v),
$$
where
$$
B^+_k(\v):=\{\sup_{k\le v\le k+1}\frac{1}{v}\x(v)\ge \v\},
~~~
B^-_k(\v):=\{\inf_{k\le v\le k+1}\frac{1}{v}\x(v)\le -\v\},
$$
we have
\begin{equation}{\label{6.4}}
{\bf P}(\x_T\not \in A_T(U,\v))\le  \sum_{k\ge UT}{\bf P}(B^+_k(\v))+
\sum_{k\ge UT}{\bf P}(B^-_k(\v)),
\end{equation}
and it is enough to bound probabilities
$$
  {\bf P}(B^+_k(\v)),~~~{\bf P}(B^-_k(\v)).
$$
We bound ${\bf P}(B^+_k(\v))$.
Since
$$
 B^+_k(\v))\subset \{\frac{1}{k}\x(k)\ge \frac{\v}{2}\}\cup
 \{\frac{1}{k}\sup_{0\le u\le 1}(\x(k+u)-\x(k))\ge \frac{\v}{2}\},
$$
we have
\begin{equation}{\label{6.5}}
{\bf P}(B^+_k(\v))\le {\bf P}(\frac{1}{k}\x(k)\ge \frac{\v}{2})+
{\bf P}(\frac{1}{k}\overline{\x}\ge \frac{\v}{2}),
\end{equation}
where $\overline{\x}:=\sup_{0\le u\le 1}\x(u)$.

We start with the first term in (\ref{6.5}). For the rate function of $\x$
$$
  \Lambda(\a):=\sup_{\l}\{\l\a-\ln {\bf E}e^{\l\x}\},~~~\a\in \R,
$$
using  exponential Chebyshev's (Chernof's) inequality
\begin{equation}{\label{6.6}}
{\bf P}(\frac{1}{k}\x(k)\ge \frac{\v}{2})\le e^{-k\Lambda(\frac{\v}{2})},
\end{equation}
noting that since ${\bf E}\x=0$),  $\Lambda(\frac{\v}{2})>0$
for any $\v>0$.

For the next term ${\bf P}(\frac{1}{k}\overline{\x}\ge \frac{\v}{2})$ in \eqref{6.5},
first verify that $\overline{\x}$
satisfies Cramer's condition $[{\bf C}_0]$.
To see this, consider a compound Poisson process $\x^+(t)$, which jumps at the same times as $\x(t)$ but the size of jumps is absolute value of the original jumps. Then for all  $t\ge 0$
$$
  \sup_{0\le v\le t}\x(v)\le \x^+(t),
$$
and
\begin{equation}{\label{6.7}}
 0\le  \overline{\x}\le  \x^+(1).
\end{equation}
It is clear that   $\x^+(1)$ satisfies $[{\bf C}_0]$, and by  (\ref{6.7})
so does  $\overline{\x}$.

Hence just as above for all
$k$ such that $\frac{\v}{2}k\ge {\bf E}\overline{\x}$
\begin{equation}{\label{6.8}}
{\bf P}(\frac{1}{k}\overline{\x}\ge \frac{\v}{2})\le e^{-\Lambda_{\overline{\x}}(\frac{\v}{2}k)}.
\end{equation}
Due to   $[{\bf C}_0]$ for
$\overline{\x}$, the function $\Lambda_{\overline{\x}}(\frac{\v}{2}k)$  grows as $k\to \infty$
faster than some linear fucntion, and we obtain from
(\ref{6.5}), (\ref{6.6}), (\ref{6.8}) that for some
$c>0$, $C<\infty$ è and all $k\ge 1$
\begin{equation}{\label{6.9}}
{\bf P}(B^+_k(\v))\le Ce^{-kc}.
\end{equation}

Obviously, a similar bound holds for  ${\bf P}(B^-_k(\v))$
and
(\ref{6.4})  yields
$$
{\bf P}(\x_T\not \in A_T(U,\v))\le  2\sum_{k\ge [UT]}Ce^{-kc}\le 2\frac{C}{1-e^{-c}}e^{-[UT]c}.
$$
 (\ref{6.1}) now follows, and
statement $(1)$ of Lemma 6.1 is proved. Statement
$(2)$ has a similar proof.
\end{proof}

In what follows we fix
$N<\infty$ and  $\v\in (0,~\frac{1}{10})$,
 constants $U\ge 2$ and  $V\ge 2$
so that (\ref{6.1}), (\ref{6.2})
hold.

For the next auxiliary result we need further notations.
Denote for a function $g\in \V$
\begin{equation}{\label{6.10}}
\overline{g}=\overline{g}(t):= \begin{cases}
         g(t),~~~t\le 2U;\\
       g(2U+),~~~t> 2U,
 \end{cases}
 \end{equation}
and for a set $B\subset \V$
$$
  \overline{B}:=\{\overline{g}:~~g\in B\}=\cup_{g\in B}\{\overline{g}\}.
$$
Then $\overline{\V}$
consists of functions  $g\in \V$, which are constant on the half-line  $(2U,\infty)$.

{\bf Lemma 6.2.}
\begin{enumerate}
\item  If $f\in A(U,\v)\cap B(2U_1,V)$  and  $2U_1\ge U$,
then $f\in B(\infty,V+\v)$.
\item  If $g\in (f)_\v$, $f\in A(U,\d)$,
then $g\in A(U+\v,\v+\d)$.
\item If $g\in (f)_\v$, $f\in B(2U,V)$, then
$g\in B(2U-\v,V+\v)$.
\item Let $g\in (f)_{10\v}$,
$f\in A(U+\v,2\v)\subset A(U+\v,10\v)$.
Then by (2) we have $g\in A(U+11\v,20\v)$.
\item If $f\in A(U+\v,2\v)$, then
\begin{equation}{\label{6.11}}
   \r(\overline{f},f)\le 4\v<5\v.
  \end{equation}
\item Let $h\in (\overline{g})_{\B,\d}$, $\overline{g}\in B(\infty,V_1)$,
then $\r(h,\overline{g})\le \d(1+V_1)$.
\item If $h\in B(\infty,V)$ and $\r(h,g)<\v$,
then $g\in B(\infty,V+\v)$.
\end{enumerate}

\bigskip

\begin{proof}Lemma 6.2. $(1)$ is obvious.

$(2)$  By  definition of $\rho$, for any  $t$
there is $(u,\hat{\b})\in \Gamma_{\hat{f}}$, such that
\begin{equation}{\label{6.12}}
  |t-u|<\v,~~~|g(t)-\hat{\b}|<\v.
\end{equation}
Therefore, for  $t\ge U+\v$ we have $u\ge U$, $|\hat{\b}|\le \d$,
$$
  |\hat{g}(t)|\le |\hat{g}(t)-\hat{\b}|+|\hat{\b}|\le \v+\d.
$$
In other words, $g\in A(U+\v,\v+\d)$.

$(3)$  By
(\ref{6.12})  if $t\le 2U-\v$ then $u\le 2U$, $|\hat{\b}|\le V$,
$$
  |\hat{g}(t)|\le |\hat{g}(t)-\hat{\b}|+|\hat{\b}|\le \v+V.
$$
so that $g\in B(2U-\v,V+\v)$.

$(4)$ is straight forward.

$(5)$ For
$t\le 2U$ we have $|\hat{\overline{f}}(t)-\hat{f}(t)|=0$,
and for $t>2U$ we have
$$|\hat{\overline{f}}(t)-\hat{f}(t)|\le
|\hat{{f}}(2U+)|+|\hat{{f}}(t)|\le 2\v+2\v=4\v.
$$
therefore
$$
  \sup_{t}|\hat{\overline{f}}(t)-\hat{f}(t)|\le 4\v,
$$
and (\ref{6.11}) follows.

$(6)$ For any $(t,\a)\in \Gamma_{h}$
there is $(u,\b)\in \Gamma_{\overline{g}}$ such that
$$
  |t-u|<\d,~~~|\a-\b|<\d.
$$
Then
$$
  |\hat{\a}-\hat{\b}|=|\frac{\a}{1+t}-\frac{\b}{1+u}|\le
  |\frac{\a}{1+t}-\frac{\b}{1+t}|+\frac{|\b|}{1+u}
  |\frac{1}{1+t}-\frac{1}{1+u}|(1+u) \le
$$
$$
   \frac{\d}{1+t} + V_1\frac{|t-u|}{1+t}< \d+\d V_1=\d(1+V_1).
$$
We have shown that for any $(t,\hat{\a})\in \Gamma_{\hat{h}}$
there is $(u,\hat{\b})\in \Gamma_{\hat{\overline{g}}}$
such that
\begin{equation}{\label{6.13}}
  |t-u| < \d(1+V_1),~~~|\hat{\a}-\hat{\b}|<\d(1+V_1).
\end{equation}
In the same way, for any  $(u,\hat{\b})\in \Gamma_{\hat{\overline{g}}}$
there is $(t,\hat{\a})\in \Gamma_{\hat{h}}$
such that (\ref{6.13}) holds.

$(7)$ is  obvious.

\end{proof}

\section {\bf Proof of Lemma 4.1.}
\begin{proof} If $J(f)=\infty$, then (\ref{4.1})  holds.
Let $J(f)<\infty$.
 Consider first an absolutely continuous $f$, for which
 $$
   J(f)=\int_0^\infty\Lambda(f'(t))dt<\infty.
 $$
 By property $II$ of Lemma 3.1 there is $U_0<\infty$
such that
 $$
   \sup_{t\ge U_0}|\hat{f}(t)|<\frac{\v}{4}.
 $$
  Denote for  $U\ge U_0$
 $$
 C(U):=\{g\in \V:~\sup_{t\le U}|\hat{g}(t)-\hat{f}(t)|<\frac{\v}{4}\},
 $$
 $$
D(U):=\{g\in \V:~\sup_{t\ge U}|\hat{g}(t)|<\frac{\v}{4}\},
 $$
and use
\begin{equation}{\label{7.1}}
(f)_\v \supset C(U)\cap D(U),
\end{equation}
which follows from
$$
 (f)_\v \supset (f)_{\hat{\r},\v}=\{g\in \V:~\hat{\r}(f,g)<\v\},
$$
which in turn follows from (\ref{2.2}) and an obvious
$$
  (f)_{\hat{\r},\v} \supset C(U)\cap D(U).
$$
We obtain from (\ref{7.1})
\begin{equation}{\label{7.2}}
{\bf P}(\x_T\in (f)_\v)\ge {\bf P}(\x_T\in C(U))-{\bf P}(\x_T\not\in D(U)).
\end{equation}
By $(1)$ of Lemma 6.1 there is $U_1\ge U_0$  such that
\begin{equation}{\label{7.3}}
\lms_{T\to \infty}\frac{1}{T}\log {\bf P}(\x_T\not \in D(U_1))\le -2J(f).
\end{equation}
On the other hand, by the local LDP on compacts which holds in the subclass of continuous functions, eg. Theorem 4.9.3 of \cite{BIB2} we have
\begin{equation}{\label{7.4}}
\lmi_{T\to \infty}\frac{1}{T}\log {\bf P}(\x_T \in C(U_1))\ge
-\int_0^{U_1}\Lambda(f'(t))dt\ge -J(f).
\end{equation}
It follows from (\ref{7.3}) and  (\ref{7.4}) that the second term in (\ref{7.2}) is $o$-little of the first.
The desired lower bound (\ref{4.1}) now follows.

If  $f\in \V$   is not absolutely continuous, then by
property $III$  of Lemma 3.1  for an arbitrary   $\d>0$ take an absolutely continuous $g$  such that
$$
  \r(g,f)<\frac{\v}{2},~~~|J(f)-J(g)|<\d.
$$
Then
$$
  {\bf P}(\x_T\in (f)_\v)\ge {\bf P}(\x_T\in (g)_{\frac{\v}{2}}),
$$
and applying the lower bound proved above, we obtain
\begin{equation}{\label{7.5}}
\lmi_{T\to \infty}\frac{1}{T}\ln {\bf P}(\x_T \in (f)_\v)\ge
\lmi_{T\to \infty}\frac{1}{T}\ln {\bf P}(\x_T \in (g)_{\frac{\v}{2}})\ge
-J(g) \ge -J(f)-\d.
\end{equation}
(\ref{4.1}) now follows and
Lemma 4.1 is proved.
\end{proof}
\section{\bf Proof of Lemma 4.2.}
\begin{proof}
For any measurable $B\subset \V$ the following clearly holds (recall that operations $\overline{f}$, $\overline{B}$ depend on $U$, see definition (\ref{6.10}))
$$
  \{\x_T\in B\}\subset \{\overline{\x}_T\in \overline{B}\},
  ~~~{\bf P}(\x_T\in B)\le {\bf P}(\overline{\x}_T\in \overline{B}).
$$
Taking $B=(f)_\v$, we obtain
\begin{equation}{\label{8.1}}
  \{\x_T\in (f)_\v\}\subset \{\overline{\x}_T\in \overline{(f)_\v}\},
  ~~~{\bf P}(\x_T\in (f)_\v)\le {\bf P}(\overline{\x}_T\in \overline{(f)_\v}).
\end{equation}

By the ELDP for processes on a compact proven in \cite{BIB4},  it follows that for any measurable $\overline{B}\subset \overline{\V}$ and any $\d>0$
\begin{equation}{\label{8.2}}
  \lms_{T\to \infty}\frac{1}{T}\ln {\bf P}(\overline{\x}_T\in \overline{B})
  \le-J(\overline{(\overline{B})_{\B,\d}}).
\end{equation}

It is easy to see (with notations from Section  6) that
$$
  {\bf P}(\x_T\in (f)_\v)\le
  {\bf P}(\x_T\not\in A(U,\v))+{\bf P}(\x_T\not\in B(2U,V)) +
$$
$$
 {\bf P}(\x_T\in (f)_\v,~~~\x_T\in A(U,\v),~~~\x_T\in B(2U,V))=:P_1+P_2+P_3,
$$
so that
\begin{equation}{\label{8.3}}
{\bf P}(\x_T\in (f)_\v)\le P_1+P_2+P_3.
\end{equation}
By $(1)$ and $(2)$ of Lemma  6.1
  $P_1$, $P_2$ admit the exponential bound  for suitable $U=U_{\v,N}$ and
   $V=V_{\v,N}$
\begin{equation}{\label{8.4}}
  P_1+P_2\le O(e^{-T(N+o(1))})~~~\mbox{as}~~~T\to \infty.
\end{equation}

We bound $P_3$.  When it is not equal to zero
by $(2)$, $(3)$ of Lemma 6.2
we have
\begin{equation}{\label{8.5}}
  f\in A(U+\v,2\v)\cap B(2U-\v,V+\v),
\end{equation}
and
$$
P_3\le P(\x_T\in (f)_\v).
$$
Using  (\ref{8.1}) we obtain
$$
  P_3\le {\bf P}(\x_T\in (f)_\v)\le
  {\bf P}(\overline{\x}_T\in \overline{(f)_\v}).
$$
Therefore by (\ref{8.2}) for any $\d>0$
\begin{equation}{\label{8.6}}
  \lms_{T\to \infty}\frac{1}{T}\ln P_3 \le -\overline{J},
  \end{equation}
where
$$
\overline{J}:=J(\overline{(\overline{(f)_\v})_{\B,\d}}).
$$
We bound
$\overline{J}$ from below,
taking into account (\ref{8.5}).

By $(5)$ of Lemma 6.2 we have
$$
 \r(f,\overline{f})<5\v,~~~(f)_\v \subset(\overline{f})_{6\v},~~~
  \overline{(\overline{(f)_\v})_{\B,\d}}\subset
  \overline{(\overline{(\overline{f})_{6\v}})_{\B,\d}},
 $$
  therefore
  $$
   \overline{J}= J(\overline{(\overline{(f)_\v})_{\B,\d}})\ge
   J(\overline{(\overline{(\overline{f})_{6\v}})_{\B,\d}}),
    $$
 and we need a lower bound for
  $$
 \overline{J}_1:= J(\overline{(\overline{(\overline{f})_{6\v}})_{\B,\d}}).
  $$
 By  (\ref{8.5}) and (\ref{6.11}),
 taking into account  $(2)$ and $(3)$ of Lemma 6.2,
 we obtain
$$
\overline{f}\in A(U+6\v,7\v)\cap B(2U-6\v,V+6\v).
$$
Note that this implies by $(1)$ of Lemma 6.2
that
\begin{equation}{\label{8.7}}
 \overline{f}\in  B(\infty,V+13\v).
\end{equation}

 Further, let $\overline{g}\in \overline{(\overline{f})_{6\v}}$, and (\ref{8.7}) to hold.
Then by  $(7)$ of Lemma 6.2, we have $\overline{g}\in B(\infty,V+19\v)$.
By $(6)$ of Lemma 6.2,
for $h\in (\overline{g})_{\B,\d}$ it holds that
$\r(h, \overline{g})<\d(1+V+19\v)$.
This means that
$$
  (\overline{(\overline{f})_{6\v}})_{\B,\d}
  \subset (\overline{(\overline{f})_{6\v}})_{\d(1+V+19\v)},~~~
 \overline{ (\overline{(\overline{f})_{6\v}})_{\B,\d}}
  \subset \overline{(\overline{(\overline{f})_{6\v}})_{\d(1+V+19\v)}}.
$$
Thus, taking $\d=\frac{\v}{(1+V+19\v)}$,
we obtain
$$
   \overline{J}\ge \overline{J}_1 =
    J(\overline{(\overline{(\overline{f})_{6\v}})_{\B,\d}})
  \ge
  J(\overline{(\overline{(\overline{f})_{6\v}})_\v}).
$$
Now we need to bound below
$$
 \overline{J}_2:=J(\overline{(\overline{(\overline{f})_{6\v}})_{\v}}).
$$

We use the following result, which will be shown later:
{\it  for any $\g>0$, $\n>0$ it holds that}
\begin{equation}{\label{8.8}}
\overline {(\overline{(\overline{f})_\g})_\n}\subset
\overline{(\overline{f})_{\g+\n}}.
\end{equation}
By (\ref{8.8})
$$
\overline{(\overline{(\overline{f})_{6\v}})_{\v}}\subset
\overline{(\overline{f})_{7\v}},~~~
\overline{J}\ge \overline{J}_2\ge J(\overline{(\overline{(\overline{f})_{6\v}})_{\v}})\ge
J(\overline{(\overline{f})_{7\v}}),
$$
and we need to bound
$$
\overline{J}_3:=J(\overline{(\overline{f})_{7\v}}).
$$

 Taking into account that
  $\r(\overline{f},f)\le 4\v$  (see $V$ Lemma 6.2),
we obtain
$$
  \overline{(\overline{f})_{7\v}}\subset\overline{(f)_{11\v}}.
$$
Therefore
$$
\overline{J}\ge \overline{J}_3 \ge
J(\overline{(f)_{11\v}}),
$$
and we need to bound
$$
\overline{J}_4:=J(\overline{(f)_{11\v}}).
$$

 Let $g\in (f)_{11\v}$,
$f\in A(U+\v,2\v)\subset A(U+\v,10\v)$.
then by  $II$ Lemma 6.2,  $g\in A(U+12\v,21\v)$.
Further, by  $V$ Lemma  6.2,
 if
$g\in A(U+12\v,21\v)\subset A(U+12\v,24\v)$
then  $\r(\overline{g},g)\le 48 \v$.
Since also $\r(g,f)<11\v$, we obtain $\r(\overline{g},f)<59\v$.
We proved that
$$
  \overline{(f)_{11\v}}\subset (f)_{59\v},~~~
 \overline{J}_4=J(\overline{(f)_{11\v}})\ge J((f)_{59\v}),
$$
so that
\begin{equation}{\label{8.9}}
  \overline{J}\ge J((f)_{59\v}).
\end{equation}
Taking into account (\ref{8.3}), (\ref{8.4}), (\ref{8.6}), (\ref{8.9}),
we obtain for any  $N<\infty$
$$
  \lms_{T\to \infty}\frac{1}{T}\ln {\bf P}(\x_T\in (f)_\v)\le
  -\min\{N,J((f)_{59\v})\}.
$$
Since $N$ is arbitrary, the Lemma is proved.

It remains to show (\ref{8.8}).
To this end, note that in the space  $\overline{\V}\subset \V$
there is the triangular inequality, for any $\overline{f}, \overline{g}, \overline{h}\in \overline{\V}$
$$
  \r(\overline{f}, \overline{g}) \le \r(\overline{f}, \overline{h})+
  \r(\overline{h}, \overline{g}).
$$
 (\ref{8.8}) follows by the triangular inequality, and the proof of Lemma 4.2 is complete.
\end{proof}

\section{\bf Proof of Lemma 4.3.}
We continue to use notations and results of the previous Sections. Fix an  $\v\in (0,0.1)$ and $N<\infty$.
By Lemma 6.1  there are  $U=U_{\v,N}<\infty$ and $V=V_{\v,N}<\infty$
so that  (\ref{6.1}) and (\ref{6.2}) hold.

\begin{proof}
 We use a result from  \cite{BIB4} that there are finitely many
$\{\overline{f}_1,\cdots,\overline{f}_M\}$
   functions$\overline{f}_i\in \overline{\V}$
such that
\begin{equation}{\label{9.1}}
\lms_{T\to \infty}\frac{1}{T}\ln{\bf P}(\x_T\not \in E)\le -N,
\end{equation}
where
$$
 E:=\cup_{i=1}^M E_i,~~~E_i:=\{g\in \V:~\r_\B(\overline{g},\overline{f}_i)<
 \frac{\v}{V}\},~~~i\in \{1,\cdots,M\}.
$$
Define $g_i:=\overline{f}_i$ for $i=1,\cdots,M$;
$$
 F:=\cup_{i=1}^MF_i,~~~F_i:=\{g\in \V:~\r(g,g_i)<5\v\}=(g_i)_{5\v}.
$$
Then
\begin{eqnarray*}
  {\bf P}(\x_T\not \in F)&\le&
  {\bf P}(\x_T\not \in A(U,\v))+
    {\bf P}(\x_T\in A(U,\v),~~~\x_T\not \in B(2U,V))+\\
   &&
  {\bf P}(\x_T\in A(U,\v),~~~\x_T \in B(2U,V),~~~\x_T\not \in E)+\\
  &&
  {\bf P}(\x_T\in A(U,\v),~~~\x_T \in B(2U,V),~~~\x_T\in E,~~~
  \x_T\not \in F)\\
  &\le&
  {\bf P}(\x_T\not \in A(U,\v))+
    {\bf P}(\x_T\not \in B(2U,V))+
  {\bf P}(\x_T\not \in E)+\\
  &&
  {\bf P}(\x_T\in A(U,\v)\cap B(2U,V)\cap E,~~~
  \x_T\not \in F)=: P_1+P_1+P_3+P_4.
  \end{eqnarray*}
Using inequalities (\ref{6.1}), (\ref{6.2}), (\ref{9.1})
we have for $j=1,2,3$
\begin{equation}{\label{9.2}}
\lms_{T\to \infty}\frac{1}{T}\log P_j\le -N.
\end{equation}

We bound $P_4$, and show that it is nil.
To this end, use the following result, which will be proven later.\\
 {\it  For any $i\in \{1,\cdots,M\}$}
\begin{equation}{\label{9.3}}
 A(U,\v)\cap B(2U,V)\cap E_i\subset F_i.
\end{equation}
Then
$$
  \{\x_T\in A(U,\v)\cap B(2U,V)\cap E\}\subset\{\x_T\in F\};
$$
and therefore
$$
  \{\x_T\in A(U,\v)\cap B(2U,V)\cap E\}\cap\{\x_T\not \in F\}\subset
  \{\x_T\in F\}\cap \{\x_T\not \in F\}=\emptyset.
$$
Thus we established that (\ref{9.3}) implies
$$
  P_4={\bf P}(\x_T\in A(U,\v)\cap B(2U,V)\cap E,~~\x_T\not \in F) =0.
$$
 Therefore bounds in (\ref{9.2})  imply the Lemma.
It remains to show (\ref{9.3}).

 Let $f\in  A(U,\v)\cap B(2U,V)\cap E_i$,
$0\le \v<0.1$, $U\ge 2$.
We show first that for all  $t\ge 2U-2$, $u\ge 2U-1$  we have
\begin{equation}{\label{9.4}}
  \frac{|f(t)|}{1+t}\le \v,~~~
  \frac{|\overline{f}(t)|}{1+t}\le \v,~~~\frac{|\overline{f}_i(u)|}{1+u}\le 3\v.
\end{equation}
The first inequality in  (\ref{9.4}) follows from
$f\in  A(U,\v)$, the second follows from the first and the definition of $\overline{f}$.
We show the third inequality in   (\ref{9.4}).

Due to  $f\in   E_i$
we have
$\r_\B(\overline{f},\overline{f}_i)<\frac{\v}{V}$,
therefore for any  $(u,\overline{f}_i(u))$
there is  $(t,\a)\in \Gamma_{\overline{f}}$ such that
$$
  |u-t|<\frac{\v}{V}, ~~~|\overline{f}_i(u)-\a|< \frac{\v}{V}.
  $$
 Therefore for  $u\ge 2U-1$ we have
 $$
    \frac{|\overline{f}_i(u)|}{1+u}\le  \frac{|\overline{f}_i(u)-\a|}{1+u}+
    \frac{|\a|}{1+u}\le
 $$
 $$
  \frac{|\overline{f}_i(u)-\a|}{1+u}+ \frac{|\a|}{1+t} +
  \frac{|\a|}{1+t}|\frac{1}{1+t}-\frac{1}{1+u}|(1+t)=
 $$
 $$
  \frac{|\overline{f}_i(u)-\a|}{1+u}+ \frac{|\a|}{1+t} +
  \frac{|\a|}{1+t}\frac{|t-u|}{(1+t)(1+u)}(1+t)\le
 $$
 $$
   \frac{\v}{V(1+u)} +\v+\v \frac{\v}{V(1+u)}\le 3\v.
 $$
   (\ref{9.4}) is proved.

 Two results follow from   (\ref{9.4}).

 $(1)_+$. {\it For any
 $(t,\hat{\a})\in \Gamma_{\hat{f}}$ for $t\ge 2U-1$
there is  $(u,\hat{\b})=(t,\hat{\b})\in \Gamma_{\hat{\overline{f}}_i}$
 such that}
 \begin{equation}{\label{9.5}}
  |t-u|<5\v,~~~|\hat{\a}-\hat{\b}|<5\v.
\end{equation}

 $(2)_+$. {\it For any
 $(u,\hat{\b})\in \Gamma_{\hat{\overline{f}}_i}$ for $u\ge 2U-1$
 there is  $(t,\hat{\a})=(u,\hat{\a})\in \Gamma_{\hat{f}}$ such that (\ref{9.5}) holds}.

Further, since $f\in B(2U,V)\cap E_i$, for any
$(t,\a)\in \Gamma_{f}$ for $t\le 2U-1$,  evidently
$(t,\a)\in \Gamma_{\overline{f}}$, therefore for any $(t,\a)\in \Gamma_{f}$ for  $t\le 2U-1$
there is  $(u,\b)\in \Gamma_{\overline{f}_i}$ such that
$$
  |t-u|<\frac{\v}{V},~~~|\a-\b|<\frac{\v}{V}.
$$
Hence
$$
 |\hat{\a}-\hat{\b}|=|\frac{\a}{1+t}-\frac{\b}{1+u}|\le
 |\frac{\a}{1+u}-\frac{\b}{1+u}|+
  |\frac{\a}{1+u}-\frac{\a}{1+t}| \le
$$
$$
  \frac{\v}{V(1+u)}+
  \frac{|\a|}{1+t}|\frac{|t-u|}{(1+t)(1+u)}|(1+t) \le
  \v+V\frac{\v}{V(1+u)}<2\v<5\v.
$$
Thus we showed the following.

$(1)_-$. {\it for any
$(t,\hat{\a})\in \Gamma_{\hat{f}}$ for $t\le 2U-1$ there is
$(u,\hat{\b})\in \Gamma_{\hat{\overline{f}}_i}$ such that  (\ref{9.5}) holds}.

Similarly, since  $f\in B(2U,V)\cap E_i$ for any
$(u,\b)\in \Gamma_{\overline{f}_i}$ for $u\le 2U-1$
there is   $(t,\a)\in \Gamma_{\overline{f}}$ such that
$$
  |t-u|<\frac{\v}{V},~~~|\a-\b|<\frac{\v}{V}.
$$
Evidently  $(t,\a)\in \Gamma_f$.
Therefore for any
$(u,\b)\in \Gamma_{\overline{f}_i}$ for $u\le 2U-1$
there is $(t,\a)\in \Gamma_{f}$ such that

$$
 |\hat{\a}-\hat{\b}|= |\frac{\a}{1+t}-\frac{\b}{1+u}|\le
 |\frac{\a}{1+u}-\frac{\b}{1+u}|+
  |\frac{\a}{1+u}-\frac{\a}{1+t}| \le
$$
$$
  \frac{\v}{V(1+u)}+\frac{|\a|}{1+t}|\frac{|t-u|}{(1+t)(1+u)}|(1+t) \le
  \v+V\frac{\v}{V(1+u)}<2\v<5\v.
$$
Hence we established

$(2)_-$. {\it For any
$(u,\hat{\b})\in \Gamma_{\hat{\overline{f}}_i}$
for $u\le 2U-1$
there is  $(t,\hat{\a})\in \Gamma_{\hat{f}}$
such that (\ref{9.5}) holds}.

It follows from  $(1)_+$, $(1)_-$; $(2)_+$, $(2)_-$
that  $\r(f,\overline{f}_i)=\r(f,g_i)<5\v$, so that
 $f\in (g_i)_{5\v}$. Consequently
(\ref{9.3}) and Lemma 4.3 are proved.

\end{proof}

\section{\bf Appendix.}

{\bf Lemma 10.1.}   Let  $ {\bf E}\x(1)=a$,
 $
  B_1:=\{f\in \V:~\sup_{0\le t\le 1}f(t)\ge 1\},
 $
and  $g_v$ be continuous piece-wise linear with
$g_v(t)=\frac{t}{v}$ for $t\in [0,v]$, and  for $t\in [v,1]$
 $g'_v(t)= a$  ($g_v(t)$ has speed $a$).
Then
\begin{equation}{\label{10.1}}
  J(B_1)=\inf_{0<v\le 1}I(g_v)=\inf_{0<v\le 1}v\Lambda(\frac{1}{v})=:v_0\Lambda(\frac{1}{v_0}).
\end{equation}

\begin{proof} Take an arbitrary $g\in B_1$ and let  $v=\e_g:=\inf\{t\in [0,1]: g(t)\ge 1\}$ be  the first  time of hitting or jumping over  level 1 and $b=\chi_g:=g(v)-1$ be     the overshoot.
Consider the function
$$
  \overline{g}(t):=g(t)-b \mathbf{1}_{(v,1]}(t),
$$
Then  $\overline{g}\in B_1$ with  the first hitting time of 1 being $v$,  but with   zero overshoot. It is clear from the definition of the rate function \eqref{J0U} that
$J(g)\ge J(\overline{g})$. Therefore we conclude that
\begin{equation}{\label{10.2}}
J(B_1)=\inf_{g\in B_1: \chi_g=0}J(g).
\end{equation}

Consider next  $g\in B_1$, with
$\e_g=v$, $\chi_g=0$, and also consider
$g_v=g_v(t)$,  defined above  in Lemma 10.1. Write the rate function as
$$
 J(g)=J_0^v(g)+J_v^1(g),\;\;\mbox{where}\;\;J_v^1(g):=J_0^1(g)-J_0^v(g),
$$
It follows from the definition of the rate function $J_0^v(g)$ (see \eqref{Thm421},Theorem 4.2.1 of \cite{BIB2}) that there exists a sequence of piece-wise linear $g_{(n)}$ on $[0,v]$ such that $g_n(v)=g(v)$ and
 $$
  J_0^v(g)=\lim_{n\to \infty}I_0^v(g_{(n)}),
$$
where
$$
  I_0^v(g_{(n)}):=\int_0^v\Lambda(g'_{(n)}(t))dt.
$$
By   convexity of $\Lambda(\a )$ and definition of the function $g_v$ for any $n$
$$ I_0^v(g_{(n)})\ge  I_0^v(g_v).$$
Taking limit as $n\to\infty$ we obtain
$$
 J_0^v(g)\ge I_0^v(g_v).
$$
Notice next that due to $\Lambda(a)=0$,
$$
  I_v^1(g_v):=I_0^1(g_v)-I_0^v(g_v)=
  \int_v^1\Lambda(g'_v(t))dt=(1-v)\Lambda(a)=0.
$$
Thus it follows
\begin{equation}{\label{10.3}}
J_0^v(g)\ge J_0^v(g_v)=I_0^v(g_v),~~~J_v^1(g)\ge J_v^1(g_v)=I_v^1(g_v).
\end{equation}
Therefore for any
  $g\in B_1$  with  $\chi_g=0$,
\begin{equation}{\label{10.4}}
  J(g)\ge J(g_v).
\end{equation}
The result now follows by  (\ref{10.2}) and (\ref{10.4}). The second equality in (\ref{10.1}) is due to
$$
  J(g_v)=I(g_v)=v\Lambda(\frac{1}{v}),
$$
so that
$$
  J(B_1)=\inf_{0<v\le 1}v\Lambda(\frac{1}{v}).
$$
Convexity of
 $
   v\Lambda(\frac{1}{v})
 $ follows from  convexity of $\Lambda$.
Hence there is
$v_0\in [0,1]$ in
(\ref{10.1}), and the proof is complete.
\end{proof}

\end{document}